\documentclass[10pt,reqno]{article}
\usepackage{latexsym, amsmath, amssymb, a4, epsfig}

\newcommand{\eqnref}[1]{(\ref {#1})}

\newcommand{\ran}{\rangle}
\newcommand{\Om}{\Omega}

\newcommand{\ds}{\displaystyle}
\newcommand{\p}{\partial}
\newcommand{\pf}{\medskip \noindent {\sl Proof}. \ }
\newcommand{\qed}{\hfill $\Box$ \medskip}
\newcommand{\RR}{\mathbb{R}}
\newcommand{\nonum}{\nonumber \\}
\newcommand{\Scal}{\mathcal{S}}
\newcommand{\Kcal}{\mathcal{K}}

\newcommand{\Lcal}{\mathcal{L}}
\newcommand{\tlambda}{\widetilde{\lambda}}
\newcommand{\tmu}{\widetilde{\mu}}
\newcommand{\tnu}{\widetilde{\nu}}
\newcommand{\tScal}{\widetilde{\mathcal{S}}}
\newcommand{\vp}{\varphi}
\def\nm{\noalign{\medskip}}
\newcommand{\pd}[2]{\frac {\p #1}{\p #2}}
\newcommand{\beq}{\begin{equation}}
\newcommand{\eeq}{\end{equation}}
\newcommand{\bfm}[1]{\mbox{\boldmath ${#1}$}}
\newcommand{\Real}{\mathop{\rm Re}\nolimits}
\newcommand{\Imag}{\mathop{\rm Im}\nolimits}

\newcommand{\Ga}{\alpha}

\newcommand{\Gd}{\delta}

\newcommand{\GT}{\Theta}

\newcommand{\GU}{\Upsilon}

\newcommand{\BGve}{\bfm\varepsilon}

\def\BS{{\bf S}}

\newtheorem{thm}{Theorem}[section]
\newtheorem{cor}[thm]{Corollary}
\newtheorem{lem}[thm]{Lemma}

\numberwithin{equation}{section}

\begin{document}
\title{Solutions to the Conjectures of P\'olya-Szeg\"o and Eshelby}

\author{Hyeonbae Kang\thanks{Department of Mathematical Sciences and RIM, Seoul National University,
Seoul 151-747, Korea}
\\ hkang@math.snu.ac.kr \and Graeme W. Milton\thanks{Department of
Mathematics, University of Utah, Salt Lake City, UT 84112, USA.}
\\ milton@math.utah.edu } \maketitle

\begin{abstract}
Eshelby showed that if an inclusion is of elliptic or ellipsoidal
shape then for any uniform elastic loading the field inside the
inclusion is uniform. He then conjectured that the converse is
true, {\it i.e.} that if the field inside an inclusion is uniform
for all uniform loadings, then the inclusion is of elliptic or
ellipsoidal shape. We call this the weak Eshelby conjecture. In
this paper we prove this conjecture in three dimensions. In two
dimensions, a stronger conjecture, which we call the strong
Eshelby conjecture, has been proved: If the field inside an
inclusion is uniform for a single uniform loading, then the
inclusion is of elliptic shape. We give an alternative proof of
Eshelby's conjecture in two dimensions using a hodographic
transformation. As a consequence of the weak Eshelby's conjecture,
we prove in two and three dimensions a conjecture of P\'olya and
Szeg\"o on the isoperimetric inequalities for the polarization
tensors. The P\'olya-Szeg\"o conjecture asserts that the inclusion
whose electrical polarization tensor has the minimal trace takes
the shape of a disk or a ball.
\end{abstract}

\bigskip



\section {Introduction}
It is well known that amongst all inclusions occupying a given
volume the sphere is the unique inclusion with maximum surface
area. This raises the question as to whether the sphere is
uniquely optimal with respect to other properties, such as
electrical properties. It was conjectured by P\'olya and Szeg\"o
\cite{PS51} that the sphere would be the unique inclusion
minimizing the trace of the electrical polarization tensor when
the inclusion and matrix have isotropic electrical properties.
Here we prove this conjecture.

A closely related conjecture is the Eshelby conjecture. It is
connected to two problems. In the transformation problem a region
(the ``inclusion'') in a homogeneous medium undergoes a
temperature change or phase change which in the absence of the
confining surrounding medium (the ``matrix'') would lead to a
uniform strain $\BGve_0$. In the elastic polarization problem, the
inclusion has different moduli to that of the matrix, and a
uniform stress is applied at infinity. If one of these problems
has been solved and the field in the inclusion is uniform then (by
subtraction or addition of a constant field) one immediately has a
solution to the other problem. Eshelby \cite{esh57,esh61},
following earlier work in special cases by Mindlin and Cooper
\cite{MC} and Robinson \cite{Rob}, showed that the stress in the
inclusion was uniform for ellipsoids and furthermore stated
(without proof) that ``among closed surfaces the ellipsoid alone
has this convenient property''. Here, for an isotropic matrix, we
prove Eshelby's conjecture (more precisely what we call the ``weak
Eshelby conjecture'') that the inclusion is necessarily
ellipsoidal if the field in the inclusion is uniform for all
transformation strains in the transformation problem, or
equivalently if the field in the inclusion is uniform for all
uniform loadings. In the transformation problem the strain
$\BGve(x)$ inside the inclusion depends linearly on $\BGve_0$ so
that we may write $\BGve(x)=\BS(x)\BGve_0$ where $\BS(x)$ is the
fourth order Eshelby tensor field. The weak Eshelby conjecture
states that if $\BS(x)$ is constant inside the inclusion then the
inclusion is an ellipsoid. It follows from certain ``trace
properties'' of the second derivative of the Green's function
associated with the problem (see, for example, equation (6.30) in
\cite{MK}) that the isotropic part of $\BS(x)$ is always uniform
and independent of the shape of the inclusion \cite{ZZD}.
Consequently for any inclusion with a sufficiently high degree of
symmetry the value of $\BS(x)$ at its center and the average of
$\BS(x)$ over the inclusion equal the value of $\BS$ in a sphere
(or circle in two-dimensions) \cite{NT97,kawa,fran}.

For planar elasticity Sendeckyj \cite{sen70}, and for antiplane
elasticity Ru and Schiavone \cite{rs96}, proved a stronger
conjecture (what we call the ``strong Eshelby conjecture'') that
the inclusion is necessarily elliptical (ellipsoidal) if the field
in the inclusion is uniform for a single transformation strain in
the transformation problem, or equivalently for a single uniform
loading in the elastic polarization problem.

Eshelby's conjecture drew increased attention when it was claimed
(see \cite{mura97,mura00} and references therein) that the field
was uniform inside star-shaped polygonal inclusions in
contradiction to the proof of Sendeckyj. Rodin \cite{rodin} proved
directly that the field cannot be uniform inside polygons or
polyhedra, and exact expressions for these non-uniform fields were
later obtained \cite{NT97,kawa,NT}. Markenscoff showed that the
field cannot be uniform if any portion of the boundary was planar
\cite{mark97} and that the only small perturbations of any
ellipsoid boundary that preserve field uniformity in the interior
are those which perturb the ellipsoid into another ellipsoid
\cite{mark98}. Lubarda and Markenscoff \cite{lubmark} showed that
the field cannot be uniform for inclusions bounded by polynomial
surfaces of higher than second degree, nor for inclusions bounded
by segments of two or more different surfaces, and argued that
non-convex inclusions are also excluded.

We remark, in passing, that not only is the field in an ellipsoid
uniform for uniform loadings but it is also polynomial for
polynomial loadings. This was proved for an isotropic matrix by
Eshelby \cite{esh61} and for an anisotropic matrix independently
by Willis in an unpublished essay \cite{willis}, and by Asaro and
Barnett \cite{AB}.

Let us now put these conjectures in a precise mathematical
framework. Consider in $\RR^d$, $d=2,3$ an inclusion $\Om$, which
is a bounded Lipschitz domain being inserted into a homogeneous
medium of conductivity $1$ in which there existed a uniform
electric field $E=-a$. We assume that the conductivity of $\Om$ is
$k \neq 1$. The insertion of the inclusion perturbs the uniform
electric field and the perturbed electric field is given by
$E=-\nabla u$ where the potential $u$ is the solution to
 \begin{equation} \label{ud}
 \ \left \{
 \begin{array}{ll}
 \ds \nabla  \cdot  \big(1 + (k-1) \chi(\Om) \big) \nabla u=0 \quad
 & \mbox{in\/ } \RR^d  ,  \\
 \nm \ds u(x)-a \cdot x=O(|x|^{1-d})  \quad & \mbox{as } |x| \to
 \infty,
 \end{array}
 \right .
 \end{equation}
where $a$ is a constant vector in $\RR^d$ indicating the direction
of the uniform field and $\chi(\Om)$ denotes the indicator
function of $\Om$. The solution $u$ to \eqnref{ud} has a multipole
asymptotic expansion at infinity, with the leading term being the
dipolar one:
 \beq \label{dipole}
 u(x)= a \cdot x + \frac{1}{\omega_d} \frac{\langle a , M x \ran}{|x|^d} +
 O(|x|^{-d}), \quad\mbox{as } |x| \to \infty.
 \eeq
Here $\omega_d$ is the area of the $d-1$ dimensional unit sphere
and $M$ is a constant $d\times d$ matrix independent of $a$ and
$x$. The matrix $M=M(\Om):=(M_{ij})$ is called the polarization
tensor associated with the inclusion $\Om$. See \cite{AKlec,
milton}.

In their book \cite{PS51} P\'olya and Szeg\"o conjectured that the
inclusion whose polarization tensor (PT) has the minimal trace
take the shape of a disk or a ball. The purpose of this paper is
to prove this conjecture in two and three dimensions. In fact, we
prove a theorem much stronger than the P\'olya and Szeg\"o
conjecture.

In connection with the P\'olya and Szeg\"o conjecture various
kinds of isoperimetric inequalities for the PT have been obtained.
See, for example, \cite{payne1, PP86, SS49}. The optimal
isoperimetric inequalities for the PT have been obtained by Lipton
\cite{Lipton93}, and later by Capdeboscq-Vogelius \cite{CV2} based
on the variational argument in \cite{kohnmilton}. The bounds are
called the Hashin-Shtrikman bounds after names of the scientists
who first found the optimal bounds on the effective conductivity
of isotropic two-phase composites \cite{hs63}, since as pointed
out in the caption of Figure 2 of \cite{MMM}, the PT bounds for
isotropic $M$ can be obtained as the low volume fraction limit of
their bounds; more generally for non-isotropic $M$ the PT bounds
can be obtained as the low volume fraction limit of the bounds of
Lurie and Cherkaev \cite{LC81,LC86} and Murat and Tartar
\cite{MT}. The PT bounds are given as follows: Let $|\Om|$ denote
the volume of $\Om$. Then
 \begin{equation}\label{HSbound-1}
 \mbox{Tr}(M)\leq |\Om| (k-1) (d-1 + \frac{1}{k}),
 \end{equation}
 and
 \begin{equation}\label{HSbound-2}
 |\Om| \mbox{Tr}(M^{-1})\leq\frac{d-1+k}{k-1},
 \end{equation}
where Tr denotes the trace.

\begin{figure}[t!]
\begin{center}
\epsfig{figure=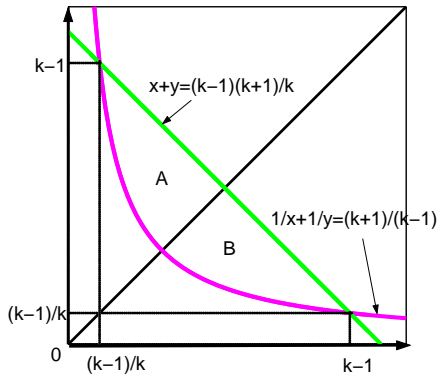, height=5cm}
\end{center}
\caption{The Hashin-Shtrikman bounds for the PTs in
$\RR^2$.}\label{Figure1}
\end{figure}

In this paper we prove the following theorem.
\begin{thm} \label{PSthm}
Let $\Om$ be a simply connected bounded Lipschitz domain in
$\RR^d$, $d=2,3$. If the polarization tensor $M(\Om)$ of $\Om$
satisfies the equality in \eqnref{HSbound-2}, then $\Om$ must be
an ellipse or an ellipsoid.
\end{thm}

Observe that if a PT $M$ has a minimal trace, then $M$ attains the
equality in \eqnref{HSbound-2} and
 \beq \label{defM}
 M= \frac{d(k-1)}{k+d-1} I
\eeq where $I$ is the $d \times d$ identity matrix assuming that
the volume $|\Om|=1$. In fact, it can be seen clearly from Figure
\ref{Figure1}, which is taken from \cite{ackkl}. In that figure,
the horizontal and vertical axis represent the eigenvalues of the
PT in two dimensions, and hence the constant trace lines are those
with slope $-1$. Thus the minimal trace occurs at the unique
tangent point of the lower hyperbola and a line with slope $-1$.
This point is an eigenvalue pair of the PT associated with the
disk. The same argument works for three dimensional case as well.
Therefore, as an immediate consequence of Theorem \ref{PSthm} we
obtain the following corollary.

\begin{cor}[The P\'olya and Szeg\"o conjecture] \label{maincor}
Let $\Om$ be a simply connected bounded Lipschitz domain in
$\RR^d$, $d=2,3$. If
 \beq \label{mintr}
 \mbox{Tr} M(\Om)= \min_{D} \mbox{Tr} M(D),
 \eeq
where $M(D)$ is the polarization tensor for the domain $D$ and the
minimum is taken over all the domain with Lipschitz boundary
(simply connected or not) with the same volume as $\Om$, then
$\Om$ is a disk or a ball.
\end{cor}

The concept of the polarization tensor appears in various contexts
such as the theory of composites (see \cite{milton} and references
therein) and the study of potential flow \cite{PS51}. Another
important usage of the concept is for the inverse boundary value
problem to detect diametrically small inclusions by means of
boundary measurements. In fact, one can approximately detect, by
boundary measurements, the location and the polarization tensor of
the inclusion. Since the polarization tensor carries important
geometric information, such as the volume of the inclusion, we are
able to recover that information from boundary measurements. It
was Friedman and Vogelius \cite{FV89} who first used the
polarization tensor for the detection of small inclusions. We
refer to \cite{AKlec} and references therein for recent
developments of this theory. It is worthwhile mentioning that the
method works for detection of multiple closely spaced inclusions
\cite{AKKL03}.

The main step in proving Theorem \ref{PSthm} is the following
theorem from \cite{km}.

\begin{thm} \label{pse}
Let $\Om$ be a bounded Lipschitz domain in $\RR^d$, $d=2,3$. If
the polarization tensor $M(\Om)$ of $\Om$ satisfies the equality
in \eqnref{HSbound-2}, then for any vector $a\in \RR^d$ the
solution $u$ to \eqnref{ud} is linear in $\Om$.
\end{thm}

Thanks to Theorem \ref{pse}, Theorem \ref{PSthm} is now an
immediate consequence of the following theorem.

\begin{thm} \label{weakeshelby}
Let $\Om$ be a simply connected bounded Lipschitz domain in
$\RR^d$, $d=2,3$. The solution $u$ to \eqnref{ud} is linear in
$\Om$ for any vector $a$ if and only if $\Om$ is an ellipse or
ellipsoid.
\end{thm}

It should be noted that only the three dimensional case in Theorem
\ref{weakeshelby} is new. In two dimensions Ru and Schiavone
\cite{rs96} proved a stronger theorem using conformal mappings: If
the gradient of the solution to \eqnref{ud} is constant for a
single non-zero direction $a$, then $\Om$ is an ellipse. This is
the anti-plane elasticity case of the strong Eshelby conjecture
for elasticity, which we now explain in the context of two and
three dimensional elasticity. Consider an elastic inclusion $\Om$,
whose Lam\'e parameters are $\tlambda, \tmu$, embedded in a medium
in $\RR^d$ with Lam\'e parameters $\lambda, \mu$. In \cite{esh57},
Eshelby showed that if $\Om$ is an ellipse or an ellipsoid, then
for any given uniform loading the elastic field inside $\Om$ is
uniform, and in \cite{esh61} he conjectured that ellipses and
ellipsoids are the only domains with this property, which is
called {\it Eshelby's uniformity property}.

In order to explain Eshelby's conjecture more precisely, let
$C=(C_{ijkl})$ be the elasticity tensor of the inclusion-matrix
composite, namely,
 \beq \label{cijkl}
 C_{ijkl} := \Bigr ( \lambda \, \chi(\RR^d \setminus \overline{\Om}) +
 \tlambda \, \chi (\Om) \Bigr ) \delta_{ij} \delta_{kl}
 + \Bigr ( \mu \, \chi(\RR^d \setminus \overline{\Om}) + \tmu\,
 \chi (\Om) \Bigr ) (\delta_{ik} \delta_{jl} + \delta_{il}
 \delta_{jk} )\;.
 \eeq
It is always assumed that
 \beq \label{con1}
 \mu > 0, \quad d \lambda + 2 \mu >0, \quad \tmu > 0
 \quad\mbox{and}\quad d \tlambda + 2 \tmu >0\;,
 \eeq
and for technical reasons we also assume that
 \beq \label{ltlmtm}
 (\lambda - \tlambda) (\mu - \tmu ) \geq 0,
 \eeq
which means $\lambda - \tlambda$ and $\mu - \tmu$ have the same
signs. For given constants $d \times d$ matrix $A$, consider the
following problem for the Lam\'e system of the linear elasticity:
 \begin{equation} \label{ud-elas}
 \ \left \{
 \begin{array}{ll}
 \ds \nabla  \cdot  \big( C (\nabla{\mathbf{u}}+\nabla{\mathbf{u}}^{T})) \big) =0 \quad
 & \mbox{in\/ } \RR^d  ,  \\
 \nm \ds \mathbf{u}(x)-Ax =O(|x|^{1-d})  \quad & \mbox{as } |x| \to
 \infty .
 \end{array}
 \right .
 \end{equation}
If $\mathbf{u}$ is the solution to \eqnref{ud-elas}, then $\nabla
\mathbf{u}$ represents the field perturbed due to the presence of
the inclusion $\Om$ under the uniform loading given by $\nabla
(Ax)$. The conductivity model \eqnref{ud} in two dimensions can be
regarded as the anti-plane elasticity model of \eqnref{ud-elas}.
What we call the strong Eshelby conjecture asserts that if the
solution $\mathbf{u}$ to \eqnref{ud-elas} for a single nonzero $A$
is linear inside $\Om$, then $\Om$ is an ellipse or an ellipsoid.
What we call the weak Eshelby conjecture states that if the
solution $\mathbf{u}$ to \eqnref{ud-elas} is linear inside $\Om$
for all $A$, then $\Om$ is an ellipse or an ellipsoid. Theorem
\ref{weakeshelby} can be regarded as a solution to the weak
Eshelby conjecture for the conductivity model. We prove the weak
Eshelby conjecture for elasticity. We only state the theorem in
three dimensions:

\begin{thm}[Weak Eshelby's conjecture in 3D] \label{weakeshelby2}
Let $\Om$ be a simply connected bounded Lipschitz domain in
$\RR^3$. The solution $\mathbf{u}$ to \eqnref{ud-elas} is linear
in $\Om$ for all $A$ if and only if $\Om$ is an ellipsoid.
\end{thm}

The strong Eshelby conjecture in two dimensions was proved by
Sendeckyj for elasticity \cite{sen70}. In this paper we give a
proof of the Eshelby conjecture in two dimensions which is
completely different from that in \cite{sen70}. The novelty of our
proof is the use of the hodographic transformation. The same
approach enables us to construct multiple inclusions satisfying
Eshelby's uniformity property \cite{kkm}.

\begin{thm}[Strong Eshelby conjecture in 2D] \label{mainthm-1}
Suppose that $d=2$. Let $\Om$ be a simply connected bounded domain
with the Lipschtz boundary. If the solution $\mathbf{u}$ to
\eqnref{ud-elas} is linear inside $\Om$ for a single nonzero $A$,
then $\Om$ must be an ellipse.
\end{thm}

We also give an alternative proof of Ru and Schiavone's theorem
for the conductivity model:
\begin{thm} \label{mainthm}
Suppose that $d=2$. Let $\Om$ be a simply connected bounded
Lipschitz domain. If the solution $u$ to \eqnref{ud} for a single
vector $a \neq 0$ is linear in $\Om$, then $\Om$ must be an
ellipse.
\end{thm}

Our proof uses a hodographic transformation. These have been
widely used to solve free boundary problems in various problems in
mechanics and fluid dynamics, to name one, the Saffman-Taylor
fingering problem \cite{st58, saf86, bp99}. It is also appropriate
to mention the Vidergauz microstructure. Vigdergauz considered a
periodic array of inclusions occupying a given volume fraction and
found the inclusion shape with minimal overall elastic energy
\cite{vig86, vig94} under certain loadings. He used the fact that
the shape would be optimal for one of these loadings if the field
inside the inclusion was uniform and hydrostatic. See \cite{gk95}
for a somewhat simpler treatment. Grabovsky and Kohn proved in the
latter paper that the low volume fraction limit of the Vigdergauz
microstructure is an ellipse. Thus one can expect that some
variant of Vigdergauz's complex analytic method might lead us to
the proof of the strong Eshelby conjecture in two-dimensions. We
regard the hodographic transformation as such a variant (see also
section 23.9 of \cite{milton}). Incidentally, we remark that for a
dilute periodic array of holes under shear loadings the ellipse is
not the optimal energy minimizing shape \cite{CGMS}.

It should be emphasized that the conjectures of P\'olya-Szeg\"o
and Eshelby are true only for simply connected domains. In a
forthcoming paper \cite{kkm}, we construct a family of structures
with two inclusions in which fields are uniform, and their PT
satisfies the lower equality in \eqnref{HSbound-2}.

This paper is organized as follows. In section 2, we review basic
facts about single layer potentials for harmonic equations and for
linear elasticity. In section 3, we prove Theorem \ref{PSthm}. We
then prove Theorem \ref{weakeshelby2} in section 3, and Theorem
\ref{mainthm-1} and \ref{mainthm} in section 4.

\section{Single layer potentials}

We review some basic facts about single layer potentials for the
harmonic equation and for isotropic elasticity. We will consider
them only in three dimensions. For details of the materials
presented here, we refer readers to \cite{AKlec}.

The single layer potential for the harmonic equation on a bounded
Lipschitz domain $\Om$ in $\RR^3$ is defined to be
 \beq
 \Scal_\Om [\phi] (x) := \frac{1}{4\pi} \int_{\p \Om} \frac{\phi(y)}{|x-y|}
 d\sigma(y), \quad x \in \RR^3,
 \label{scal}
 \eeq
where $\phi$ is a square integrable function on $\p\Om$ and
$d\sigma(y)$ is the surface measure. Thus $\Scal_\Om [\phi]$ is a
function on $\RR^3$ and $\Scal_\Om [\phi](x)$ denotes its value at
$x$. The following boundary behavior of the normal derivative of
the single layer potential is well-known:
 \beq  \label{singlejumpcond}
 \pd{}{n} \Scal_\Om [\phi] \bigg |_\pm (x) = \biggl(\pm \frac{1}{2}
 I + \Kcal_\Om^*\biggr) [\phi] (x) \quad \mbox{ a.e. } x \in \partial
 \Om ,
 \eeq
where $n=(n_1,n_2,n_3)$ is the outward unit normal to $\Om$,
$\pd{}{n}$ denotes the normal derivative, and $\Kcal_\Om^*$ is
defined by
 \beq
 \Kcal^* _\Om [\phi] (x) = \frac{1}{4\pi} \mbox{p.v.} \int_{\p \Om}
 \frac{\langle x-y, n(x) \rangle}{|x-y|^3} \phi(y) \, d \sigma(y).
 \eeq
Here the subscripts $+$ and $-$ denote the limits from the outside
and inside $\Om$, respectively, and p.v. denotes the Cauchy
principal value. See \cite{Folland76} for a proof of
\eqnref{singlejumpcond} when $\p \Om$ is smooth and \cite{ver84}
when $\p \Om$ is Lipschitz. It is known \cite{KS2000} (see also
\cite[Section 2.47]{AKlec}) that the solution $u$ to \eqnref{ud}
is given by
 \beq \label{ujx}
 u(x) = a \cdot x + \Scal_\Om [\phi](x), \quad x \in \RR^3,
 \eeq
where
 \beq
 \phi = \left( \frac{k+1}{2(k-1)} I- \Kcal^*_\Om \right)^{-1}[a \cdot n] \quad \mbox{on } \p
 \Om .
 \eeq
Furthermore, we have
 \beq \label{phij}
 \phi = (k-1) \pd{u}{n} \Big|_{-}.
 \eeq
The invertibility of the operator $\frac{k+1}{2(k-1)} I-
\Kcal^*_\Om$ on $L^2(\p \Om)$ is established in \cite{EFV92}.

It is worthwhile to note that in view of the jump  relation
\eqnref{singlejumpcond}, the usage of the single layer potential
is natural since \eqnref{ud} when $d=3$ is equivalent to the
following problem:
 \begin{equation} \label{ud-2}
 \ \left \{
 \begin{array}{ll}
 \ds \Delta u=0  \quad
 & \mbox{in\/ } \Om \cup (\RR^3 \setminus \overline{\Om})  ,  \\
 \nm \ds u|_+=u|_- \quad & \mbox{on\/ } \p \Om, \\
 \nm \ds \pd{u}{n} \Big|_+= k \pd{u}{n} \Big |_- \quad & \mbox{on\/ } \p \Om, \\
 \nm \ds u(x)-a \cdot x=O(|x|^{-2})  \quad & \mbox{as } |x| \to
 \infty,
 \end{array}
 \right .
 \end{equation}

We now review a similar representation formula for isotropic
elasticity. The elastostatic system corresponding to the Lam{\'e}
constants $\lambda, \mu$ is defined by
 \beq \label{defmulam}
 \Lcal_{\lambda, \mu} \mathbf{u} := \mu \Delta \mathbf{u} + (\lambda +
 \mu) \nabla (\nabla \cdot \mathbf{u})\;.
 \eeq
The corresponding conormal derivative ${\partial
\mathbf{u}}/{\partial \nu}$ on $\p \Om$ is defined to be
\begin{equation} \label{Lame-2}
\pd{\mathbf{u}}{\nu} := \lambda (\nabla \cdot \mathbf{u}) n + \mu
(\nabla \mathbf{u} + \nabla \mathbf{u}^T ) n \quad \mbox{on } \p
\Om \;,
\end{equation}
where the superscript $T$ denotes the transpose of a matrix. The
Kelvin matrix ${\bf \Gamma} = ( \Gamma_{ij} )_{i,j=1}^d$ of the
fundamental solution  to the Lam{\'e} system $\Lcal_{\lambda,
\mu}$ in three dimensions is given by
 \beq \label{Kelvin}
 \Gamma_{ij} (x) : = - \frac{\Ga_1}{4 \pi} \frac{\delta_{ij}}{|x|} -  \frac{\Ga_2}{4 \pi}
 \frac{x_i x_j}{|x|^3}, \quad  x \neq 0\;,
 \eeq
where
 \beq \label{ab} \Ga_1= \frac{1}{2} \left ( \frac{1}{\mu} +
 \frac{1}{2\mu + \lambda} \right ) \quad\mbox{and}\quad \Ga_2=
 \frac{1}{2} \left ( \frac{1}{\mu} - \frac{1}{ 2 \mu + \lambda}
 \right )\;.
 \eeq
The single layer potentials of the density function $\vp$ on $\p
\Om$ associated with the Lam{\'e} parameters $(\lambda, \mu)$ are
defined by
 \beq
 \vec{\Scal}_\Om [\mathbf{\vp}] (x)  := \int_{\p \Om} {\bf \Gamma} (x-y)
 \mathbf{\vp} (y)\, d \sigma (y)\;, \quad x \in \RR^3\;.
 \label{singlelayer}
 \eeq
The single layer potential enjoys the following jump relation:
 \beq \label{singlejumplame}
 \pd{}{\nu} \vec{\Scal}_\Om [\mathbf{\vp}] \Big|_+ - \pd{}{\nu} \vec{\Scal}_\Om
 [\mathbf{\vp}]
 \Big|_- = \vp \quad \mbox{on } \p \Om,
 \eeq
where $\partial /{\partial \nu}$ denotes the conormal derivative
defined in \eqnref{Lame-2}.

Let $\Psi$ be the vector space of all linear solutions of the
equation $\Lcal_{ \lambda , \mu } \mathbf{u} =0$ and ${\partial
\mathbf{u}}/{\partial \nu} =0$ on $\p \Om$, or alternatively, \beq
\label{defpsi}
 \Psi := \bigg\{ \mathbf{\psi} : \p_i \psi_j + \p_j \psi_i =0, \quad
 1 \leq i, j \leq d \bigg\}\;.
\eeq Here the $\psi_i$ for $i=1,\ldots,d,$ denote the components
of $\mathbf{\psi}$. Define
 \beq \label{l2psi}
 L^2_\Psi (\p \Om) := \bigg \{ ~ \mathbf{f} \in L^2 (\p \Om) : ~ \int_{\p \Om} \mathbf{f} \cdot
 \mathbf{\psi}\,
 d\sigma =0 \ \mbox{for all } \mathbf{\psi} \in \Psi ~ \bigg\}
 \eeq
which is a subspace of codimension $6$ in $L^2(\partial \Om)$.
Note that \eqnref{ud-elas} is equivalent to the following problem:
 \begin{equation} \label{trans-6}
 \begin{cases}
 \ds \Lcal_{ \lambda , \mu } \mathbf{u} = 0
 \quad \mbox{in } \RR^3 \setminus \overline{\Om}\;, \\
 \nm  \ds \Lcal_{\tlambda , \tmu} \mathbf{u} = 0 \quad \mbox{in } \Om\;,  \\
 \nm  \ds
 \mathbf{u} \big |_{-} = \mathbf{u} \big |_{+} \quad \mbox{on } \p \Om\;\;, \\ \nm  \ds
 \pd{\mathbf{u}}{\tnu} \bigg |_{-} = \pd{\mathbf{u}}{\nu} \bigg |_{+}
 \quad \mbox{on } \p \Om\;, \\ \nm
 \ds \mathbf{u} - Ax = O(|x|^{-2}) \quad \mbox{as } |x| \to \infty,
 \end{cases}
 \end{equation}
where $\Lcal_{\tlambda , \tmu}$ and $\tnu$ are the Lam\'e operator
and the conormal derivative with respect to the Lam\'e constants
$(\tlambda, \tmu)$ of the inclusion.  We denote by
$\vec{\Scal}_\Om$ and $\vec{\tScal}_\Om$ the single layer
potentials on $\p \Om$ corresponding to the Lam{\'e} constants
$(\lambda, \mu)$ and $(\tlambda, \tmu)$ of the matrix and
inclusion, respectively. We then have the following representation
formula for the solution $\mathbf{u}$ to \eqnref{ud-elas} or
equivalently \eqnref{trans-6}: There exists a unique pair
$(\mathbf{\vp}, \mathbf{\psi})  \in L^2(\p \Om) \times L^2_{\Psi}
(\p \Om)$ such that the solution $\mathbf{u}$ of \eqnref{trans-6}
is represented by
 \begin{equation} \label{rep-1}
 \mathbf{u} (x) =
 \begin{cases}
 Ax + \vec{\Scal}_\Om [\mathbf{\psi}] (x)\;, \quad & x \in \RR^3 \setminus
 \overline \Om\;, \\ \vec{\tScal}_\Om [\mathbf{\vp}] (x)\;, \quad & x \in \Om\;,
 \end{cases}
 \end{equation}
where the pair $(\mathbf{\vp}, \mathbf{\psi})$  is the unique
solution in $L^2 (\p \Om) \times  L^2_{\Psi} (\p \Om)$ of
 \beq \label{ES93-2}
 \begin{cases}
 \ds \vec{\tScal}_\Om [\mathbf{\vp}] \big |_{-} - \vec{\Scal}_\Om [\mathbf{\psi}] \big |_{+} =
 (Ax)|_{\p \Om} \quad & \mbox{on } \p \Om\;, \\
 \nm \ds \pd{}{\tnu} \vec{\tScal}_\Om [\mathbf{\vp}] \bigg|_{-} - \pd{}{\nu}
 \vec{\Scal}_\Om [\mathbf{\psi}] \bigg |_{+} = \pd{(Ax)}{\nu} \bigg |_{\p \Om} \quad & \mbox{on } \p
 \Om\;.
 \end{cases}
 \eeq
The unique solvability of the integral equation \eqnref{ES93-2}
was proved in \cite{ES93} under the condition \eqnref{ltlmtm},
which is why we assumed this condition.

\section{Proof of the P\'olya-Szeg\"o conjecture}

We now prove Theorem \ref{weakeshelby}. Theorem \ref{PSthm}
follows as an immediate consequence. We consider only the three
dimensional case because the same proof works for the two
dimensional case. We begin with the following lemma.

\begin{lem} \label{atob}
Suppose that $d=3$ and that the solution $u$ to \eqnref{ud} is
linear in $\Om$ for any applied field $a\in \RR^3$. Define a
linear transformation $\Lambda$ on $\RR^3$ by $\Lambda(a)= \nabla
u_a|_\Om$ where $u_a$ is the solution to \eqnref{ud}. Then
$\Lambda$ is one-to-one and onto.
\end{lem}

\pf Since the equation in \eqnref{ud} is linear, that $\Lambda$ is
linear is obvious. Thus it suffices to show that $\Lambda$ is
one-to-one. Suppose that $\Lambda(a)=0$ for some $a \in \RR^3$.
Then $\nabla u_a=0$ in $\Om$. Define $v(x)= u_a(x)-a\cdot x$ for
$x \in \RR^3$. Then $v$ is the solution to
 \begin{equation} \label{ud-10}
 \ \left \{
 \begin{array}{ll}
 \ds \nabla  \cdot  \big(1 + (k-1) \chi(\Om) \big) \nabla v=0 \quad
 & \mbox{in\/ } \RR^d  ,  \\
 \nm \ds v(x) =O(|x|^{1-d})  \quad & \mbox{as } |x| \to
 \infty,
 \end{array}
 \right .
 \end{equation}
By the uniqueness of the solution to \eqnref{ud-10}, $v \equiv 0$
and hence $a=0$. This completes the proof. \qed

Lemma \ref{ud-10} can be interpreted as follows: For any vector $b
\in \RR^3$ there exists $a$ such that the solution $u$ to
\eqnref{ud} satisfies
 \beq
 u(x) = b \cdot x + c, \quad x \in \Om
 \eeq
for some constant $c$. It then follows from \eqnref{ujx} and
\eqnref{phij} that
 \beq
 (b-a) \cdot x + c = (k-1) \Scal_\Om [b \cdot n] (x), \quad x \in
 \Om.
 \eeq
In other words, $\Scal_\Om [b \cdot n]$ is linear in $\Om$ for any
$b \in \RR^3$. In particular, we have
 \beq \label{somnj}
 \Scal_\Om [n_j] (x) = \mbox{linear in } \Om, \quad j=1,2,3.
 \eeq
By reversing the arguments one can see that \eqnref{somnj} is
equivalent to the solution $u$ to \eqnref{ud} being linear in
$\Om$ for any vector $a$.

Now although we have only defined the action of the functional
$\Scal_\Om$ on scalar functions, the obvious generalization of
\eqnref{scal} defines its action on vector valued functions. In
particular we have
 \beq \label{scalb}
 \Scal_\Om[n](x)= -\nabla \int_\Om \frac{1}{4\pi |x-y|} dy, \quad x
 \in \Om,
 \eeq
which can be seen using the divergence theorem. Thus we get from
\eqnref{somnj}
 \beq \label{intB}
 \int_\Om \frac{1}{4\pi |x-y|} dy = \mbox{a quadratic polynomial}, \quad x
 \in \Om .
 \eeq

Since the property \eqnref{intB} is independent of the
conductivity ratio $k \neq 1$, we have an interesting consequence.
\begin{cor}
For $j=1,2,3$ and $0< k \neq 1< \infty$, let $u_j^{(k)}$ be the
solution to \eqnref{ud} with $a=\mathbf{e}_j$ and conductivity
ratio $k$ where $\{ \mathbf{e}_1, \mathbf{e}_2, \mathbf{e}_3 \}$
is the standard basis for $\RR^3$. If $\nabla u_j^{(k)}$ is
constant in $\Om$ for $j=1,2,3$ and for some $k$, then $\nabla
u_j^{(k)}$ is constant in $\Om$ for all $k$.
\end{cor}

Suppose that \eqnref{intB} holds, {\it i.e.}, there is a symmetric
matrix $A$, a constant vector $b$, and a constant $C$ such that
\beq \label{newt}
 \frac{1}{4\pi} \int_\Om \frac{1}{|x-y|} dy = x \cdot Ax + b \cdot x+ C, \quad
 x\in \Om.
 \eeq
This identity has an interesting physical interpretation. If we
think of $\Om$ as a body of constant density in free space then
the left hand side is (to within a proportionality factor) just
the Newtonian gravitational potential at $x$, and the identity
\eqnref{newt} then says that the gradient of this potential, which
is the gravitational field, depends linearly on $x$ within the
body. It has been shown by  Dive \cite{div31} and Nikliborc
\cite{nik32} that ellipsoids are the only bodies which have this
property. Before explaining their proof, let us establish some
elementary results.

After a unitary transformation if necessary, we may assume $A$ is
diagonal and
 \beq \label{sumj13}
 \frac{1}{4\pi} \int_\Om \frac{1}{|x-y|} dy = \frac{1}{2} \sum_{j=1}^3 a_j x_j^2 + b \cdot x+ C, \quad
 x\in \Om,
 \eeq
for some constants $a_j$ and $C$, and a constant vector $b$. We
claim that each $a_j$ is positive. In fact, it follows from
\eqnref{scalb} that
 \beq \label{scalajxj}
 \Scal_\Om [n_j](x) = -a_j x_j - b_j, \quad x \in \Om, \quad j=1,2,3.
 \eeq
Therefore, by \eqnref{singlejumpcond}, we have
 \beq \label{ajnonzero}
 (-\frac{1}{2} I + \Kcal_\Om^*)[n_j] = -a_j n_j \quad \mbox{on } \p \Om,
 \eeq
and hence
 \beq
 \left( \frac{k+1}{2(k-1)} I - \Kcal_\Om^* \right)^{-1}[n_j]= \frac{k-1}{1+ (k-1) a_j} n_j \quad \mbox{on } \p \Om.
 \eeq
Since $(-\frac{1}{2} I + \Kcal_\Om^*)$ is invertible on $L^2_0(\p
\Om):= \{ f \in L^2(\p \Om) : \int_{\p \Om} f =0 \}$ as was proved
in \cite{ver84}, one can see from \eqnref{ajnonzero} that $a_j
\neq 0$. Moreover, the polarization tensor $M=(M_{ij})$ associated
with $\Om$ is given by
 \beq
 M_{ij} = \int_{\p \Om} y_j \left( \frac{k+1}{2(k-1)} I - \Kcal_\Om^* \right)^{-1}[n_i](y)
 d\sigma(y)
 \eeq
as was proved in \cite{AK01}. Therefore, we have
 \beq
 M_{ij}= \delta_{ij} |\Om| \frac{k-1}{1+(k-1) a_j}, \quad i,j=1,2,3.
 \eeq
Since the polarization $M$ is positive definite if $k>1$ and
negative definite if $k<1$ regardless of $k$ (see \cite{AKlec}),
we have $a_j >0$, $j=1,2,3$.

We can now make a complete square out of \eqnref{sumj13} and make
a translation if necessary to conclude that
 \beq \label{sumfinal}
 \frac{1}{4\pi} \int_\Om \frac{1}{|x-y|} dy = \frac{1}{2} \sum_{j=1}^3 a_j x_j^2 + C, \quad
 x\in \Om,
 \eeq
with $a_j >0$.

The following theorem was proved by Dive \cite{div31} and
Nikliborc \cite{nik32} for a $C^1$ domain. The same theorem for
Lipschitz domains can be proved by a slight variation of their
arguments.

\begin{thm} \label{dive}
Let $\Om$ be a bounded domain with a Lipschitz boundary. The
relation \eqnref{sumfinal} holds if and only if $\Om$ is an
ellipsoid of the form
 \beq \label{c1c2c3}
 \frac{x^2_1}{c^2_1}+\frac{x^2_2}{c^2_2}+\frac{x^2_3}{c^2_3} \le 1.
 \eeq
\end{thm}

\pf We briefly sketch the proof. Note that if $\Om$ is an
ellipsoid of the form \eqnref{c1c2c3}, then \eqnref{sumfinal}
holds with
 \beq \label{ajcj}
 a_j = \frac{c_1c_2c_3}{2}\int^{\infty}_{0}\frac{ds}{(c_j^2+s)\sqrt{(c_1^2+s)(c_2^2+s)(c_3^2+s)}},
 \eeq
for $j=1,2,3$ \cite{div31, nik32}.

To prove the converse,  Suppose \eqnref{sumfinal} holds with $a_j
>0$, $j=1,2,3$. Then there is a unique triple $c_1, c_2, c_3$
satisfying the relation \eqnref{ajcj}. Existence of such a triple
was proved in \cite{div31, nik32}. Let $E$ be the ellipsoid given
by \eqnref{c1c2c3}. Then, defining for any region $\GU$ \beq
\label{defN} N_\GU(x):=\frac{1}{4\pi}
 \int_\GU \frac{1}{|x-y|} dy,
\eeq we have
 \beq \label{intom}
 N_E(x)= \frac{1}{2} \sum_{j=1}^3 a_j x_j^2 +C_1,  \quad x
\in \overline{E},
 \eeq
for some constant $C_1$. For $t>0$, let $E_t:= \{ tx | x \in E
\}$. Then by simple scaling one can see that
 \beq \label{intomt}
 N_{E_t}(x) = \frac{1}{2} \sum_{j=1}^3 a_j x_j^2 +C_t,  \quad x
\in \overline{E_t},
 \eeq
for some constant $C_t$ depending only on $t$. Let $t_0$ be the
smallest number such that $\Om \subset E_t$ for all $t \ge t_0$.
Then there is a point $Q$ which is contained in $\p E_{t_0} \cap
\p \Om$, and
 \beq \label{eom}
 N_{E_{t_0} \setminus \Om}(x) = N_{E_{t_0}}(x) -
 N_\Om(x)= \mbox{constant}, \quad x \in \Om.
 \eeq

Let $n(Q)$ be the unit outward normal to $E_{t_0}$ at $Q$. Since
$E_{t_0} \setminus \Om$ lies in one side of the tangent plane to
$E_{t_0}$ at $Q$, we have
 \beq \label{nqnu}
 \nabla N_{E_{t_0} \setminus \Om}(Q) \cdot n(Q) = - \frac{1}{4\pi}
 \int_{E_{t_0} \setminus \Om} \frac{\langle Q-y, n(Q) \rangle}{|Q-y|^3} dy < 0,
 \eeq
provided that $E_{t_0} \setminus \Om$ is not empty.

If $\p\Om$ is $C^1$, then $n(Q)$ is also normal to $\Om$ at $Q$
and the normal line goes through $\Om$. Here and in what follows,
``goes through $\Om$'' means that there is $s_0>0$ such that the
line segment $\{ Q-s n(Q)| 0<s<s_0 \} \subset \Om$. But by
\eqnref{eom} we get $\nabla N_{E_{t_0} \setminus \Om}(Q) \cdot
n(Q)=0$, and hence we can conclude that $E_{t_0} \setminus \Om =
\emptyset$ by \eqnref{nqnu}. This is the argument in \cite{div31,
nik32}.

If $\p\Om$ is only Lipschitz, then we can argue as follows. By
\eqnref{nqnu}, we have $\nabla N_{E_{t_0} \setminus \Om}(Q) \neq
0$ provided that $E_{t_0} \setminus \Om$ is not empty. Thus for
any unit vector $v_0$ and for any open neighborhood $V$ of $v_0$
in $S^2$, the unit sphere, there is $v \in V$ such that
 \beq \label{nvne}
 \nabla N_{E_{t_0} \setminus \Om}(Q)\cdot v \neq 0.
 \eeq
Choose $v_0$ so that the line in the direction $v_0$ passing
through $Q$ goes through $\Om$. Since $\p\Om$ is Lipschitz, there
is a neighborhood $V$ in $S^2$ of $v_0$ such that any line in the
direction of the vector in $V$ passing through $Q$ goes through
$\Om$. Then  for some $v \in V$ \eqnref{nvne} holds. But by
\eqnref{eom} we get contradiction and hence $E_{t_0} \setminus \Om
= \emptyset$. This completes the proof of Theorem \ref{dive} and
hence Theorem \ref{weakeshelby}.

\section{Proof of the weak Eshelby conjecture}

In this section we prove Theorem \ref{weakeshelby2}, the weak
Eshelby conjecture for elasticity. We only will prove the three
dimensional case because again the same proof works for the two
dimensional case.

In fact there is a close link between the weak Eshelby conjecture
for elasticity and the weak Eshelby conjecture for conductivity
(which have just proved). In composite microstructures of two
isotropic phases it is known that if for some periodic
microgeometry the elastic field is uniform and hydrostatic (i.e.
proportional to the identity) in phase one then that
microstructure necessarily attains the ``bulk modulus type trace
bound'', equation (6.35) in \cite{MK}, or the opposite inequality,
depending on the moduli of the phases. Then as a consequence of an
argument of Grabovsky \cite{grab96} (see also section 25.6 of
\cite{milton}) a solution to the conductivity problem (for any
direction of the applied field) can be generated from that
elasticity field, and the electric field in phase one is also
necessarily uniform. This argument strongly suggests that if
Eshelby's uniformity property holds for all applied loadings, then
there will be one loading for which the field in the inclusion is
hydrostatic, and from which one can generate solutions to the
conductivity problem. As a consequence the electric field in the
inclusion will be uniform for all applied uniform electric fields,
and hence the inclusion shape must be ellipsoidal.

Let us see this directly. Suppose that $d=3$ and that the solution
$\mathbf{u}$ to \eqnref{ud-elas} is linear in $\Om$ for any $3
\times 3$ matrix $A$. Then in the representation formula
\eqnref{rep-1}, $\vec{\tScal}_\Om \mathbf{\vp}$ is linear in
$\Om$, say
 \beq
 \vec{\tScal}_\Om [\mathbf{\vp}] (x) = Bx + b, \quad x \in \Om,
 \eeq
for some $3 \times 3$ matrix $B$ and vector $b$. It then follows
from the integral equation \eqnref{ES93-2} that
 \beq \label{ES93-3}
 \begin{cases}
 \ds \vec{\Scal}_\Om [\mathbf{\psi}] \big |_{+} = Bx-Ax+b \quad & \mbox{on } \p \Om\;, \\
 \nm \ds \pd{}{\nu} \vec{\Scal}_\Om [\mathbf{\psi}] \bigg |_{+} =
 \ds \pd{(Bx)}{\tnu} \bigg|_{-} - \pd{(Ax)}{\nu} \bigg |_{\p \Om} \quad & \mbox{on } \p
 \Om\;.
 \end{cases}
 \eeq
Since $\Lcal_{\lambda, \mu} \vec{\Scal}_\Om [\mathbf{\psi}]=
\Lcal_{\lambda, \mu} (Bx-Ax+b)=0$ in $\Om$, the first relation in
\eqnref{ES93-3} implies that
 \beq
 \vec{\Scal}_\Om [\mathbf{\psi}] (x) = Bx-Ax+b, \quad x \in \Om.
 \eeq
It then follows from the jump relation \eqnref{singlejumplame}
 \beq \label{psirel}
 \psi = \pd{}{\nu} \vec{\Scal}_\Om [\mathbf{\psi}] \bigg |_{+} - \pd{}{\nu} \vec{\Scal}_\Om [\mathbf{\psi}] \bigg
 |_{-} = \pd{(Bx)}{\tnu} - \pd{(Bx)}{\nu} \quad \mbox{on } \p\Om.
 \eeq
Substituting \eqnref{psirel} into the first identity in
\eqnref{ES93-3}, we arrive at
 \beq \label{Bline}
 \vec{\Scal}_\Om \left[ \pd{(Bx)}{\tnu} - \pd{(Bx)}{\nu}
 \right] (x)= \mbox{linear}, \quad x \in \Om.
 \eeq

The following theorem can be proved in the exactly same manner as
Lemma \ref{atob}.

\begin{lem} \label{atob-elas}
Suppose that $d=3$ and that the solution $\mathbf{u}$ to
\eqnref{ud-elas} is linear in $\Om$ for any $3 \times 3$ matrix
$A$. Define a linear transformation $\Lambda$ on $V$, the vector
space of all $3 \times 3$ matrices, by $\Lambda(A)= \nabla
\mathbf{u}_A|_\Om$ where $\mathbf{u}_A$ is the solution to
\eqnref{ud-elas}. Then $\Lambda$ is one-to-one and onto.
\end{lem}

Since Lemma \ref{atob-elas} means that for any $3 \times 3$ matrix
$B$, there is $A$ such that the solution $\mathbf{u}$ to
\eqnref{ud-elas} satisfies $\nabla \mathbf{u}=B$ in $\Om$, we can
conclude that \eqnref{Bline} holds for all $3 \times 3$ matrix
$B$.

Let us take a hydrostatic field $B$, i.e. with $B_{pq}=\Gd_{pq}$
the identity matrix, in \eqnref{Bline}. Then we have $Bx= x$. We
claim that
 \beq \label{pd13}
 \vec{\Scal}_{\Om} \left[ \pd{x}{\nu} \right]
 (x)_j = - \frac{2\mu+3 \lambda}{2\mu+\lambda} \Scal_\Om [n_j](x),
 \quad x \in \Om, \ \
 j=1,2,3.
 \eeq
where the subscript $j$ indicates the $j$-th component. Similarly,
we have
 \beq \label{pd14}
 \vec{\Scal}_{\Om} \left[ \pd{x}{\nu} \right]
 (x)_j = - \frac{2\tmu+3 \tlambda}{2\mu+\lambda} \Scal_\Om [n_j](x),
 \quad x \in \Om, \ \
 j=1,2,3.
 \eeq
We will give proofs of \eqnref{pd13} and \eqnref{pd14} at the end
of this section.

Because of \eqnref{pd13} and \eqnref{pd14}, we get
 \beq \label{tmumu}
\vec{\Scal}_\Om \left[ \pd{x}{\tnu} - \pd{x}{\nu}
 \right](x)_j = \frac{(2\mu+3 \lambda)-(2\tmu+3 \tlambda)}{2\mu+\lambda} \Scal_\Om
 [n_j](x), \quad x \in \Om.
 \eeq
Notice that the constant on the right-hand side of \eqnref{tmumu}
is not zero unless $2\mu+3 \lambda-(2\tmu+3 \tlambda)$ which is
excluded by our assumption \eqnref{ltlmtm}. Since the left-hand
side of \eqnref{tmumu} is linear in $\Om$, we deduce that
$\Scal_\Om [n_j]$ is linear in $\Om$ for $j=1,2,3$. We now
conclude from the result in the previous section that $\Om$ is an
ellipsoid, and the proof is complete.

Let us now prove \eqnref{pd13}. Since
\begin{align}
\pd{x}{\nu} &= \lambda \nabla \cdot (x) n + \mu (\nabla (x) +
\nabla (x)^T ) n = (2\mu+3\lambda)n, \label{xxx}
\end{align}
where $n=(n_1,n_2,n_3)$ is the unit outward normal, we get, for
$j=1,2,3$, \beq \label{sxpeqk}
 \vec{\Scal}_{\Om} \left[ \pd{x}{\nu} \right] (x)_j
 = (2\mu+3\lambda)\left ( \int_{\p \Om} \mathbf \Gamma(x-y)n(y) \right )_j
= (2\mu+3\lambda) \sum_{\ell=1}^3 \int_{\p \Om}
\Gamma_{j\ell}(x-y) n_\ell(y). \eeq Using \eqnref{Kelvin} we have
theorem,
\begin{align}
& \sum_{\ell=1}^3 \int_{\p \Om} \Gamma_{j\ell}(x-y) n_\ell(y)
d\sigma(y)
\nonumber \\
&= -\frac{\Ga_1}{4\pi} \int_{\p \Om} \frac{n_j(y)}{|x-y|}
d\sigma(y) - \frac{\Ga_2}{4\pi} \int_{\p \Om} (x_j-y_j)
\frac{\langle x-y, n(y) \rangle}{|x-y|^{3}} d\sigma(y).
\label{410}
\end{align}
By Green's theorem, we have
\begin{align}
& \int_{\p \Om} (x_j-y_j) \left[\pd{}{n_y} \frac{1}{|x-y|}\right]
d\sigma(y)- \int_{\p \Om} \left[\pd{}{n_y}(x_j-y_j)\right]
\frac{1}{|x-y|} d\sigma(y) \nonum & = \int_{\Om} (x_j-y_j)
\left[\Delta_y \frac{1}{|x-y|}\right] dy - \int_{\Om}
\left[\Delta_y (x_j-y_j)\right] \frac{1}{|x-y|} dy =0,
\label{green}
\end{align}
and hence
 \beq \label{410-1}
 \int_{\p \Om} (x_j-y_j) \frac{\langle x-y, n(y) \rangle}{|x-y|^{3}} d\sigma(y) =
 - \int_{\p \Om} \frac{n_j(y)}{|x-y|} d\sigma(y).
 \eeq
It then follows from \eqnref{410} and the fact that  $\Ga_1-\Ga_2=
\frac{1}{2\mu+ \lambda}$ that \beq \sum_{l=1}^3 \int_{\p \Om}
\Gamma_{jl}(x-y) n_l(y) d\sigma(y) = \frac{-\Ga_1+\Ga_2}{4\pi}
\int_{\p \Om} \frac{n_j(y)}{|x-y|} d\sigma(y) = - \frac{1}{(2\mu+
\lambda)}\Scal_\Om[n_j](x) . \label{diagonal} \eeq By substituting
this back in \eqnref{sxpeqk} we get \eqnref{pd13}.

To prove \eqnref{pd14}, it suffices to note that
\begin{align}
\pd{x}{\nu} &= \tlambda \nabla \cdot (x) n + \tmu (\nabla (x) +
\nabla (x)^T ) n = (2\tmu+3\tlambda)n. \label{xxxx}
\end{align}
The proof is now complete.

\section{Alternative Proofs of the strong Eshelby conjecture in 2D}

In this section we give alternative proofs of Theorem
\ref{mainthm-1} and \ref{mainthm}. We first prove a lemma on
univalence of the analytic functions which will be used to prove
Theorem \ref{mainthm-1} and \ref{mainthm}. The proof of the
following lemma relies on the level curve argument which was used
in various contexts. We particularly mention the work of
Alessandrini and Nesi \cite{an01} in which they showed the
univalence of $\sigma$-harmonic mappings in the context of
periodic composite materials. The two phase inclusion-matrix
problem in the free space can be viewed as a low volume limit of
the periodic case.

\begin{lem} \label{unival}
Let $\Om$ be a simply connected domain with Lipschitz boundary in
$\mathbb{C}$ and let $f$ be a analytic function in $\mathbb{C}
\setminus \overline{\Om}$ such that there are constants $\alpha
\neq 0$ and $\beta$ such that $f(z)-(\alpha z +\beta) \to 0$ as
$|z| \to \infty$. If $f(z)= i x_2$ for $z=x_1+ix_2 \in \p \Om$,
then $f$ is univalent in $\mathbb{C} \setminus \overline{\Om}$.
\end{lem}

\pf Let $u$ and $v$ be the real and imaginary parts of $f$,
respectively. Observe that $f$ maps $\mathbb{C}^* \setminus
\overline{\Om}$ onto $\mathbb{C}^* \setminus [ic_1,ic_2]$ where
$c_1=\min_{z\in \p \Om} v(z)$ and $c_2=\max_{z\in \p \Om} v(z)$,
and $\mathbb{C}^*$ is the Riemann sphere. Since $\mathbb{C}^*
\setminus \overline{\Om}$ and $\mathbb{C}^* \setminus [ic_1,ic_2]$
are simply connected, it suffices to show that $f'(z) \neq 0$ for
$z \in \mathbb{C}^* \setminus \overline{\Om}$ to prove univalence
of $f$. It is obvious that $f'(\infty) \neq 0$. We will show that
$\nabla u(z) \neq 0$ for any $z \in \mathbb{C} \setminus
\overline{\Om}$.

Suppose that $\nabla u(z_0)=0$ for some $z_0 \in \mathbb{C}
\setminus \overline{\Om}$. Then there is an integer $m \ge 2$ such
that $u$ takes the form
 \beq \label{expan}
 u(z)-u(z_0)= \sum_{n=m}^\infty |z-z_0|^n (a_n \cos n \theta + b_n
 \sin n \theta), \quad z \mbox{ near } z_0,
 \eeq
where $z-z_0= |z-z_0| e^{i\theta}$ and $a_m^2+b_m^2 \neq 0$.
Therefore, there are $2m$ branches of level curves $C(z_0):=\{~ z~
|~ u(z)=u(z_0) ~\}$ coming out of the point $z_0$. Since $u(z)=
ax_1+bx_2+c+O(|z|^{-1})$ as $|z| \to \infty$ for some real
constants $a,b$, and $c$, at most two of these $2m$ branches are
unbounded and extend to the infinity. For the other branches of
$C(z_0)$ which are bounded, one of the following two occurs: (i) a
branch intersects $\p \Om$, (ii) a branch meets another branch and
makes a loop.

Suppose that (ii) occurs. Then the loop may encircle
$\overline{\Om}$ or be contained in $\mathbb{C}\setminus
\overline{\Om}$. If it is contained in $\mathbb{C}\setminus
\overline{\Om}$, then by the maximum principle, $u$ is constant
inside the loop and we have contradiction. If the loop encircles
$\overline{\Om}$, then $u$ is constant (zero) on $\p \Om$ and on
the loop (the constants may be different). Let $\GT$ be the
annular region enclosed by the loop and $\p \Om$. Then by the
Cauchy-Riemann equations, we have $\pd{v}{n}=0$ on $\p \GT$, and
hence $v$ is constant in $\GT$, which leads us to a contradiction.

Now suppose that (ii) does not occur for any bounded branch. Let
$l_j$, $j=1, \ldots, k$ ($k\ge 2$) be bounded branches of
$C(z_0)$. Let $z_j$ be the point where $l_j$ intersects $\p \Om$.
Then there are two branches, say $l_1$ and $l_2$, such that $l_1$,
$l_2$, and the arc on $\p \Om$ connecting $z_1$ and $z_2$ make the
boundary of a connected region $\Om$. Since $u= u(z_1)=u(z_2)=0$
on $l_1$ and $l_2$, it follows from the maximum principle that
$u=0$ in $\Om$, which contradicts to the assumption that $\alpha
\neq 0$. This completes the proof. \qed

We now provide the alternative proofs of Theorem \ref{mainthm-1}
and \ref{mainthm}.

\noindent {\sl Proof of Theorem \ref{mainthm}}. Let $u$ be the
solution to \eqnref{ud} for some $a \neq 0$, and suppose that $u$
is linear in $\Om$. Put $u^i:= u|_{\Om}$ and $u^e:=u|_{\RR^2
\setminus \overline{\Om}}$. It is proved in \cite{IP90} that there
are functions $U^i$ and $U^e$ analytic in $\Om$ and $\RR^2
\setminus \overline{\Om}$, respectively, such that $\Real U^i=u^i$
and $\Real U^e=u^e$, and
 \begin{equation} \label{uiue-2}
 \frac{k+1}{2} U^i - \frac{k-1}{2} \overline{U^i} = U^e + i c \quad \mbox{on } \p
 \Om,
 \end{equation}
for some real constant $c$. Since $u^i$ is linear in $\Om$, there
are complex numbers $\delta$ and $\gamma$ such that \beq
\label{Ulin}
 U^i(z) = \delta z+\gamma.
\eeq Define $f_\Om$ by
 \beq
 f_\Om(z): = U^e(z) - \frac{k+1}{2} (\delta z + \gamma) + i c, \quad
 z \in \mathbb{C}\setminus \overline{\Om}.
 \eeq
Then by \eqnref{uiue-2} we have
 \beq \label{uiue-3}
 f_\Om(z) = - \frac{k-1}{2} (\overline{\delta} \overline{z} +
 \overline{\gamma} ), \quad z \in \p \Om.
 \eeq
We now define $\psi_\Om$ by
 \begin{equation}
 \psi_\Om(z)= \frac{1}{\overline{\delta}(k-1)} \left[ f_\Om(z) + \frac{k-1}{2} \overline{\gamma} \right]
 + \frac{1}{2} z, \quad z \in \mathbb{C}\setminus \overline{\Om}.
 \end{equation}
Then one can see from \eqnref{uiue-3} that
 \begin{equation}
 \psi_\Om(z)= ix_2, \quad \mbox{for } z=x_1+ix_2 \in \p \Om.
 \end{equation}
Since $u^e(x)-a \cdot x = O(|x|^{-1})$ as $|x| \to \infty$,
$\psi_\Om$ takes the form
 \beq
 \psi_\Om(z)=\alpha z + \beta + \vp_\Om(z),
 \eeq
for some complex numbers $\alpha$ and $\beta$, and a analytic
function $\vp_\Om$ in $\mathbb{C} \setminus \overline{\Om}$ such
that $\vp_\Om(z)=O(|z|^{-1})$ as $|z| \to \infty$.

We claim that $\alpha \neq 0$. In fact, if $\alpha =0$, then $h:=
\Real \psi_\Om$ is the solution to
 \beq \label{const}
 \begin{cases}
 \Delta h = 0 \quad & \mbox{in } \RR^2 \setminus \overline{\Om}, \\
 h= 0 &\mbox{on } \p \Om, \\
 h(x) - c = O(|x|^{-1}) \quad &\mbox{as } |x| \to \infty,
 \end{cases}
 \eeq
for some real constant $c$. But the maximum principle (with the
point at infinity being regarded as a point on the Riemann sphere)
implies $c=0$ and hence $h \equiv 0$ in $\RR^2 \setminus
\overline{\Om}$. Thus we conclude that $\alpha \neq 0$.

It then follows from Lemma \ref{unival} that $\psi_\Om$ is a
univalent mapping from $\mathbb{C}^* \setminus \overline{\Om}$
onto $\mathbb{C}^* \setminus [ic_1,ic_2]$ where $c_1=\min_{z\in \p
\Om} \Imag \psi_\Om(z)$ and $c_2=\max_{z\in \p \Om} \Imag
\psi_\Om(z)$. We may assume that $c_1=-1$ and $c_2=1$ by scaling
if necessary.

If $D$ is the unit disk, we can construct such a mapping $\psi_D$
explicitly. In fact, it is the Koebe function \beq \label{koebe}
 \psi_D(z) = \frac{1}{2}(z-\frac{1}{z}), \quad |z|>1.
\eeq Then $\psi_D$ is a univalent mapping from $\mathbb{C}^*
\setminus \overline{D}$ onto $\mathbb{C}^* \setminus [-i,i]$.

Let \beq \label{defF}
 F(z) :=\psi_\Om^{-1} \circ \psi_D (z), \quad |z|>1,
\eeq in which $\psi_\Om^{-1}$ is the hodographic transformation.
Then $F$ is a univalent mapping from $|z|>1$ onto
$\mathbb{C}\setminus \overline{\Om}$ and $F(\infty)=\infty$.
Moreover, $F(z)$ behaves as a linear analytic function at
infinity. Since $\mathbb{C}^* \setminus \overline{D}$ is simply
connected, it follows from Caratheodory's theorem
\cite[P.18]{Pomm} that $F$ extends to $\mathbb{C}^* \setminus D$
as a continuous function and $F$ is a homeomorphism from $\p D$
onto $\p \Om$. Observe that if $|z|=1$, then $F(z) = \psi_\Om^{-1}
(ix_2)$, and hence $\Imag F(z) = x_2$. In other words, $\Imag F$
is a solution to
 \beq \label{imf}
 \begin{cases}
 \Delta u = 0 \quad & \mbox{in } \RR^2 \setminus \overline{D}, \\
 u(x)=x_2 &\mbox{on } \p D, \\
 u(x) - (b_1 x_1+ b_2 x_2 + c) = O(|x|^{-1}) \quad &\mbox{as } |x| \to \infty,
 \end{cases}
 \eeq
where $b_1$, $b_2$, and $c$ are real constants.

We claim that $c=0$ and
 \beq \label{imagi}
 \Imag F(z)= b_1 x_1 + b_2 x_2 + \frac{-b_1 x_1 +
 (1-b_2)x_2}{|z|^2}, \quad |z|>1.
 \eeq
In fact, if we define $u$ by
 \beq
 u(z)= \Imag F(z)- \left[ b_1 x_1 + b_2 x_2 + \frac{-b_1 x_1 +
 (1-b_2)x_2}{|z|^2} \right], \quad |z| >1,
 \eeq
then $u$ is a solution to \eqnref{const} with $\Om$ replaced by
$D$, and hence $c=0$ and $u \equiv 0$.

We now get from \eqnref{imagi} that
 \beq
 F(z)= \gamma z + C + \frac{\beta}{z}, \quad |z|>1,
 \eeq
for some constants $\gamma \neq 0$ and $\beta$. One can easily see
that the image of the unit disk under $F$ is an ellipse. This
completes the proof of Theorem \ref{mainthm}. \qed

\medskip

\noindent{\sl Proof of Theorem \ref{mainthm-1}}. Let $\mathbf{u}$
be the solution to \eqnref{ud-elas} for some constant $a_{ij}$,
not all zero, and assume that $\mathbf{u}$ is linear in $\Om$. We
first invoke the following complex representation of the solution
to \eqnref{ud-elas} from \cite{Musk77} (see also \cite[Theorem
6.20]{AKlec}): Let $\mathbf{u}=(u,v)$ be the solution of
\eqnref{ud-elas} for $d=2$ and let $\mathbf{u}_e :=
\mathbf{u}|_{\mathbb{C} \setminus \overline{\Om}}$ and
$\mathbf{u}_i := \mathbf{u}|_{\Om}$. Then there are unique
functions
 $\vp_e$ and $\psi_e$ analytic in $\mathbb{C} \setminus \overline
 {\Om}$ and $\vp_i$ and $\psi_i$ analytic in $\Om$ such that
 \begin{align}
 2\mu (u_e + i v_e)(z) & = \kappa \vp_e(z) - z \overline{\vp_e'(z)}
 - \overline{\psi_e(z)}\;,  \quad z \in \mathbb{C} \setminus \overline {\Om}\;,
 \label{ext} \\
 2\tmu (u_i + i v_i)(z) &= \widetilde\kappa \vp_i(z) - z
 \overline{\vp_i'(z)} - \overline{\psi_i(z)}\;,  \quad z \in  \Om\;,
 \label{int}
 \end{align}
 where
 \beq \label{defkappa}
  \kappa= \frac{\lambda + 3
  \mu}{\lambda+\mu}\;, \quad \quad \widetilde{\kappa} =
  \frac{\tlambda + 3 \tmu}{\tlambda+\tmu}\;.
  \eeq
 Moreover, the following holds on $\p \Om$:
 \begin{align}
 \frac{1}{2\mu} \bigg( \kappa \vp_e(z) - z \overline{\vp_e'(z)} -
 \overline{\psi_e(z)} \bigg ) & = \frac{1}{2\tmu} \bigg (
 \widetilde\kappa \vp_i(z) - z \overline{\vp_i'(z)} -
 \overline{\psi_i(z)} \bigg )\;, \label{dir} \\
 \vp_e(z) + z \overline{\vp_e'(z)} + \overline{\psi_e(z)} & =
 \vp_i(z) + z \overline{\vp_i'(z)} + \overline{\psi_i(z)} + c\;,
 \label{tra}
 \end{align}
 where $c$ is a constant. Equation \eqnref{dir} is the continuity of the
 displacement and \eqnref{tra} is the continuity of the traction.

 It follows from \eqnref{dir} and \eqnref{tra} that
  \beq \label{kaplus}
  (\kappa+1) \vp_e(z) = \Big(\frac{\mu \widetilde\kappa}{\tmu} +1
  \Big) \vp_i(z) + \Big(1-\frac{\mu}{\tmu} \Big) \Big( z \overline{\vp_i'(z)} + \overline{\psi_i(z)}
  \Big) +c,  \quad z \in \p \Om.
  \eeq
Since $\mathbf{u}$ is linear in $\Om$, $\vp_i$ and $\psi_i$ are
linear analytic functions in $\Om$ by the uniqueness of $\vp_i$
and $\psi_i$. It then follows from \eqnref{kaplus} that there are
complex numbers $\alpha, \beta, C$ such that
 \beq
 \vp_e(z) = \alpha z - \overline{\beta z} + C, \quad z \in \p \Om.
 \eeq

Suppose that $\beta \neq 0$. Define $f_\Om$ by
 \beq
 f_\Om(z) = \frac{1}{2\beta} \big[ \vp_e(z)- \alpha z - C \big ]
 + \frac{1}{2} z, \quad z \in \mathbb{C} \setminus \overline{\Om}.
 \eeq
Then $f_\Om(z)= i x_2$ for $z=x_1+ix_2 \in \p \Om$. Following the
same argument as in the proof of Theorem \ref{mainthm}, we
conclude that $\Om$ is an ellipse.

If $\beta=0$, then $\vp_e$ can be extended to the whole of
$\mathbb{C}$ as an entire function. In fact, if we define
$\vp_e(z) = \alpha z + C$ in $\Om$, then $\vp_e$ is analytic in
$\mathbb{C} \setminus \p\Om$ and continuous on $\p \Om$, and hence
is an entire function (It can be proved using Morera's theorem
\cite{ahlfors} that that if a function defined in an open set $U$
is analytic in $U$ minus a Lipschitz curve and is continuous in
$U$, then that function is analytic in $U$.) But, since
$\mathbf{u}(x)- \sum_{ij} a_{ij} x_i \mathbf{e}_j=O(|x|^{-1})$ as
$|x| \to \infty$, $\vp_e$ takes the form $\vp_e(z) =\gamma z +
f_e(z)$ for some constant $\gamma$ and a analytic function $f_e$
in $\mathbb{C} \setminus \overline{\Om}$ such that
$f_e(z)=O(|z|^{-1})$ as $|z| \to \infty$. So $f_e$ extends as a
bounded entire function and hence $f_e$ is constant and the
constant is $0$. Thus
 \beq \label{vpelin}
 \vp_e(z)=\gamma z \quad \mbox{for } z
 \in \mathbb{C} \setminus \overline{\Om}.
 \eeq
It then follows from \eqnref{dir} and \eqnref{tra} that
 \begin{align}
 \frac{\kappa-\bar \gamma}{2\mu} z  - \frac{1}{2\mu}
 \overline{\psi_e(z)} & = \frac{1}{2\tmu} \bigg (
 \widetilde\kappa \vp_i(z) - z \overline{\vp_i'(z)} -
 \overline{\psi_i(z)} \bigg )\;,  \\
 (1+\bar \gamma)z + \overline{\psi_e(z)} & =
 \vp_i(z) + z \overline{\vp_i'(z)} + \overline{\psi_i(z)} + c\;,
 \end{align}
on $\p \Om$. Therefore there are complex constant $\delta, \eta,
C$ such that
 \beq
 \psi_e(z) = \delta z - \overline{\eta z} + C.
 \eeq
If $\eta =0$, then by the same reasoning as to derive
\eqnref{vpelin}, one can see $\psi'(z)$ is constant. Therefore,
$u_e+iv_e$ is linear, and hence $\mathbf{u}(x)- \sum_{ij} a_{ij}
x_i \mathbf{e}_j=0$ for all $x \in \mathbb{C} \setminus
\overline{\Om}$. This is possible only when $\Om$ is an empty set.
Thus, $\eta \neq 0$. Now let
 \beq
 g_\Om(z) = \frac{1}{2\eta} \big[ \psi_e(z)- \delta z - C \big ]
 + \frac{1}{2} z, \quad z \in \mathbb{C} \setminus \overline{\Om}.
 \eeq
Then $g_\Om(z)= i x_2$ for $z=x_1+ix_2 \in \p \Om$. Following the
same argument as in the proof of Theorem \ref{mainthm}, we
conclude that $\Om$ is an ellipse. This completes the proof. \qed

We finally mention that the strong Eshelby conjecture for three
dimensions (except in the special case where the field in the
inclusion is hydrostatic) has not been proven, not even for the
conductivity case.

{\bf Acknowledgement}. We would like to thank Victor Isakov for
informing us of the existence of the papers \cite{div31} and
\cite{nik32}, Dave Barnett and Peter Schiavone for drawing our
attention to the papers \cite{sen70} and \cite{rs96}, Yves
Capdeboscq for stimulating discussion on the P\'olya-Szeg\"o
conjecture, and Hyundae Lee for pointing out an error in earlier
draft of this paper. H.K.is grateful for partial support by the
grant KOSEF R01-2006-000-10002-0 and G.W.M is grateful for support
from the National Science Foundation through grant DMS-0411035.

\end{document}